\documentclass[10pt]{amsart}
\usepackage{color}
\usepackage{todonotes}
\usepackage{array}
\usepackage{geometry}
\usepackage{graphicx}
\usepackage{amssymb,amsmath,amsfonts}
\usepackage{epstopdf}
\usepackage{subfig}
\usepackage[ruled,vlined,norelsize,linesnumbered]{algorithm2e}
\usepackage{ifthen,xifthen,xparse}
\DeclareDocumentCommand{\R}{ O{ } O{ } }{
	\ifthenelse{\isempty{#1}}
		{\mathbb{R} }
		{
			\ifthenelse{\isempty{#2}}
				{\mathbb{R}^{#1} }
				{ \mathbb{R}^{#1 \times #2} }
		}
}

\newcommand{\rope}{\sqrt{1+\epsilon}}

\newcommand{\vc}{\textmd{vec}}
\newcommand{\diag}{\textmd{diag}}

\theoremstyle{definition}

\title{On matrix balancing and eigenvector computation}

\author{Rodney James}
\address[Rodney James]{Colorado College, Colorado Springs, CO, USA}
\email{Rodney.James@ColoradoCollege.edu}
\author{Julien Langou$^1$}
\address[Julien Langou]{University of Colorado Denver, Denver, CO, USA}
\email{Julien.Langou@ucdenver.edu}
\author{Bradley R. Lowery$^1$}
\address[Bradley R. Lowery]{University of Colorado Denver, Denver, CO, USA}
\email{Bradley.Lowery@ucdenver.edu}

\date{\today}

\begin{document}
\maketitle

\footnotetext[1]{Research of this author was fully supported by the National Science Foundation grant \# NSF CCF 1054864.}

\begin{abstract}
	Balancing a matrix is a preprocessing step while solving the nonsymmetric eigenvalue problem.
	Balancing a matrix reduces the norm of the matrix and hopefully this will improve the accuracy of 
	the computation. 
	Experiments have shown that balancing can improve the accuracy of the computed eigenvalues.  
	However, there exists examples where balancing increases the eigenvalue condition number
	(potential loss in accuracy), deteriorates eigenvector accuracy, and deteriorates the backward
	error of the eigenvalue decomposition. 
	In this paper we propose a change to the stopping criteria of the LAPACK balancing algorithm, 
	\texttt{GEBAL}.  The new stopping criteria is better
	at determining when a matrix is nearly balanced.  Our experiments show that the new algorithm is 
	able to maintain good backward error, while improving the eigenvalue accuracy when possible.  
	We present stability analysis, numerical experiments, 
	and a case study to demonstrate the benefit of the new stopping criteria.   
\end{abstract}

\section{Introduction}

For a given vector norm $\|\cdot\|$, an $n$-by-$n$ (square) matrix $A$ is said to
be {\it balanced} if and only if, for all $i$ from 1 to $n$, the norm of its
$i$-th column and the norm of its $i$-th row are equal. 
$A$ and $\widetilde{A}$ are {\em diagonally similar} means that there exists a
diagonal nonsingular matrix $D$ such that $\widetilde{A} = D^{-1} A D $.
Computing $D$ (and/or $\widetilde{A}$) such that $\widetilde{A}$ is balanced is called {\em balancing} $A$, and
$\widetilde{A}$ is called {\em the balanced matrix}.

Using the 2-norm, given an $n$-by-$n$ (square) matrix $A$,
Osborne~\cite{Osborne:1960:PM:321043.321048} proved that there exists a unique
matrix $\widetilde{A}$ such that $A$ and $\widetilde{A}$ are diagonally similar
and $\widetilde{A}$ is balanced. Osborne~\cite{Osborne:1960:PM:321043.321048}
also provides a practical iterative algorithm for computing $\widetilde{A}$.

Parlett and Reinsch~\cite{parlett:1969:NM} (among other things) proposed to approximate $D$ with $D_2$, a
diagonal matrix containing only powers of two on the diagonal. While $D_2$ is only an approximation 
of $D$, $\widetilde{A}_2 = D_2^{-1} A D_2 $ is ``balanced enough'' for all practical matter. The key idea
is that there is no floating-point error in computing $\widetilde{A}_2$ on a base-two computing machine.

Balancing is a standard preprocessing step while solving the nonsymmetric eigenvalue
problem (see for example the subroutine
DGEEV in the LAPACK library).  Indeed, since, $A$ and $\widetilde{A}$ are
similar, the eigenvalue problem on $A$ is similar to the eigenvalue problem on
$\widetilde{A}$.  Balancing is used to hopefully improve the accuracy of the
computation~\cite{templatesSolutions:2000:siam,chen:1998:masters,chen:2001:phd,chen:2000:LAA}.
However, there are examples where the
eigenvalue accuracy~\cite{watkins:2006:case}, eigenvector
accuracy~\cite{LAPACK-forum-4270}, and backward error will deteriorate.  We
will demonstrate these examples throughout the paper.

What we know balancing does is decrease the norm of the matrix.  Currently, the stopping criteria
for the LAPACK balancing routine, \texttt{GEBAL}, is if a step in the algorithm does not provide a
significant reduction in the norm of the matrix (ignoring the diagonal elements).  Ignoring the 
diagonal elements is a reasonable assumption since they are unchanged by a diagonal similarity transformation.  
However, it is our observation that this stopping criteria
is not strict enough and allows a matrix to be ``over balanced'' at times. 

We propose a simple fix to the stopping criteria.  Including the diagonal elements and balancing with respect to 
the 2-norm solves the problem of deteriorating the backward error in our test cases.  Using the
2-norm also prevents balancing a Hessenberg matrix, which previously was an example of 
balancing deteriorating the eigenvalue condition number \cite{watkins:2006:case}.

Other work on balancing was done by Chen and Demmel~\cite{chen:1998:masters,chen:2001:phd, chen:2000:LAA}. 
They extended the balancing literature to sparse matrices.  
They also looked at using a weighted norm to balance the matrix.
For nonnegative matrices, if the weight vector is taken to be the Perron vector (which is nonnegative),
then the largest eigenvalue has perfect conditioning.  However, little can be said about the
conditioning of the other eigenvalues.  This is the only example in the literature 
that guarantees balancing improves the condition number of any of the eigenvalues.

The rest of the paper is organized as follows. In Section~\ref{sec:background} 
we provide background and review the current literature on balancing.  
In Section~\ref{sec:casestudy} we introduce a simple example
to illustrate why balancing can deteriorate the backward error (and the accuracy of the eigenvectors).
In Section~\ref{sec:backerror} we provided backward error analysis.
Finally, in Section~\ref{sec:experiments} we look at a few test cases 
and conclude that that using the 2-norm in the balancing algorithm is the
best solution.

\section{The Balancing Algorithm}\label{sec:background}

Algorithm~\ref{alg:osborne} shows Osborne's original balancing algorithm,
which balances a matrix in the 2-norm~\cite{Osborne:1960:PM:321043.321048}.  
The algorithm assumes that the matrix is {\it irreducible}.  A matrix is {\it reducible} if
there exists a permutation matrix, $P$, such that, 
\begin{equation}\label{eq:reduce}
P A P^T = \left( \begin{array}{cc} A_{11} & A_{12} \\ 0 & A_{22} \end{array} \right),
\end{equation}
where $A_{11}$ and $A_{22}$ are square matrices.  A diagonally similarity transformation can
make the off-diagonal block, $A_{12}$, arbitrarily small with $D = \diag(\alpha I, I)$ for 
arbitrarily large $\alpha$.  For converges of the balancing algorithm it is necessary for the
elements of $D$ to be bounded.  For reducible matrices, the diagonal blocks can be balanced
independently.  
\begin{algorithm}[htbp]
	\SetKwInOut{input}{Input}
	\SetKwInOut{output}{Output}
	\SetKwInOut{notes}{Notes}
	\Indm
	\input{An irreducible matrix $A \in \R[n][n]$.}
	\notes{$A$ is overwritten by $D^{-1} A D$, and converges to unique balanced matrix. $D$ also
		converges to a unique nonsingular diagonal matrix.} 
	\BlankLine
	\Indp
	$D \leftarrow I$ \\	
	\For{ $k \leftarrow 0,1,2,\dots$ } {
		\For{ $i \leftarrow 1,\dots, n$ } {
			$\displaystyle c \leftarrow \sqrt{ \sum_{j\ne i} |a_{j,i}|^2 }, \;\; 
			r \leftarrow \sqrt{ \sum_{j\ne i} |a_{i,j}|^2 }$ \\ \label{alg1.line1}
			$\displaystyle f \leftarrow \sqrt{\frac{r}{c}} $ \\ \label{alg1.line2}
			
			$d_{ii} \leftarrow f \times d_{ii}$ \\
			$ A(:,i) \leftarrow f \times A(:,i), \;\;
			A(i,:) \leftarrow A(i,:)/f$ \\
		}
	}
	\caption{ Balancing (Osborne)}
	\label{alg:osborne}
\end{algorithm}

The balancing algorithm operates on columns and rows
of $A$ in a cyclic fashion. 
In line~\ref{alg1.line1}, the 2-norm (ignoring the diagonal element) 
of the current column and row are calculated. 
Line~\ref{alg1.line2} calculates the scalar, $f$, that will be used to update $d_{ii}$,
column $i$, and row $i$.  The quantity $f^2 c^2 + r^2 / f^2$ is minimized at $f = \sqrt{\frac{r}{c}}$,
and has a minimum value of $2rc$.
It is obvious that for this choice of $f$, the Frobenius norm is non-increasing since,
 \[ \| A^{(nk+i)} \|_F^2 - \| A^{(nk+i+1)} \|_F^2 = ( c^2 + r^2 ) - ( 2rc ) = ( c - r )^2 \ge 0. \] 
Osborne showed that Algorithm~\ref{alg:osborne} 
converges to a unique balanced matrix and $D$ converges to a unique (up to a scalar multiple) nonsingular 
diagonal matrix. The balanced matrix has minimal Frobenius norm among all diagonal similarity transformations.  

Parlett and Reinsch~\cite{parlett:1969:NM} generalized the balancing algorithm to any $p$-norm.
Changing line~\ref{alg1.line1} to 
\[ 
	c \leftarrow \Big( \sum_{j\ne i} |a_{j,i}|^p \Big)^{1/p}, \;\; 
	r \leftarrow \Big( \sum_{j\ne i} |a_{i,j}|^p \Big)^{1/p},
\]
the algorithm will converge to a balanced matrix in the $p$-norm.  Parlett and Reinsch however, do not 
quantify if the norm of the matrix will be reduced or minimized.  The norm will in fact
be minimized, only it is a non-standard matrix norm.  Specifically, the balanced matrix, $\widetilde{A}$,
\[
\| \widetilde{A} \|_p = \min \left( \| \vc( D^{-1} A D ) \|_p \mid D \in \mathcal{D} \right),
\]
where $\mathcal{D}$ is the set of all nonsingular diagonal matrices and $\vc$ 
stacks the columns of $A$ into an $n^2$ column vector.  For $p = 2$, this is the Frobenius norm.

Parlett and Reinsch also
restrict the diagonal element of $D$ to be powers of the radix base (typically 2).  
This restriction ensures there is no computational error in the balancing algorithm. 
The algorithm provides only an approximately balanced matrix.  Doing
so without computational error is desirable and the exact balanced matrix is not 
necessary. Finally, they introduce a stopping criteria for the algorithm.  Algorithm~\ref{alg:parlett}
shows Parlett and Reinsch's algorithm for any $p$-norm.  Algorithm~\ref{alg:parlett} 
(balancing in the $1$-norm) is essentially what is currently implemented by LAPACK's 
\texttt{GEBAL} routine.
\begin{algorithm}[htbp]
	\SetKwInOut{input}{Input}
	\SetKwInOut{output}{Output}
	\SetKwInOut{notes}{Notes}

	\Indm
	\input{ A matrix $A \in \R[n][n]$.  We will assume that $A$ is irrreducible.}
	\output{ A diagonal matrix $D$ and $A$ which is overwritten by $D^{-1} A D$. } 
	\notes{ $\beta$ is the radix base. 
		On output $A$ is nearly balanced in the $p$-norm. }
	\BlankLine
	\Indp
	$D \leftarrow I$ \\
	$converged \leftarrow 0$ \\
	\While{ $converged = 0$ } {
		$converged \leftarrow 1$ \\
		\For{ $i \leftarrow 1,\dots,n$ } {
			$\displaystyle c \leftarrow \Big( \sum_{j\ne i} |a_{j,i}|^p \Big)^{1/p}, \;\; 
			r \leftarrow \Big( \sum_{j\ne i} |a_{i,j}|^p \Big)^{1/p}$ \\ \label{alg2.line1}
			$s \leftarrow c^p + r^p , \;\;
			f \leftarrow 1 $ \\ 			
			\While{ $c < r/\beta$ } { \label{alg2.line2}
				$\displaystyle c \leftarrow c\beta, \;\;
				r \leftarrow r/\beta, \;\;
				f \leftarrow f \times \beta $ \\
			}
			\While{ $c \ge r\beta$ } {
				$\displaystyle c \leftarrow c/\beta, \;\;
				r \leftarrow r\beta, \;\;
				f \leftarrow f/\beta $ \\ \label{alg2.line3}
			} 
			\If{ $(c^p + r^p) < 0.95 \times s$ } { \label{alg2.line4}
				$ converged \leftarrow 0, \;\; 
				d_{ii} \leftarrow f \times d_{ii}$ \\
				$ A(:,i) \leftarrow f \times A(:,i), \;\;
				A(i,:) \leftarrow A(i,:)/f$ \\
			}
		}	
	}
	\caption{ Balancing (Parlett and Reinsch) }
	\label{alg:parlett}
\end{algorithm}

Again, the algorithm proceeds in a cyclic fashion operating on a single row/column at a time.
In each step, the algorithm seeks to minimize $f^p c^p + r^p / f^p$.
The exact minimum is obtained for $f = \sqrt{r/c}$.
Lines~\ref{alg2.line2}-\ref{alg2.line3} find the closest approximation to the exact minimizer.
At the completion of the second inner while loop, $f$ satisfies the inequality
\[ \beta^{-1}\left( \frac{r}{c} \right) \le f^2 < \left( \frac{r}{c}\right)\beta. \] 
Line~\ref{alg2.line4} is the stopping criteria.  It states that the current step must decrease
$f^p c^p + r^p / f^p$ to at least 95\% of the previous value.  If this is not satisfied then the step is skipped.
The algorithm terminates when a complete cycle does not provide significant decrease.

Our observation is that the stopping criteria is not a good indication of the true decrease in the 
norm of the matrix.  If the diagonal element of the current 
row/column is sufficiently larger than the rest of the entries, then balancing will be
unable to reduce the matrix norm by a significant amount (since the diagonal element is unchanged).
Including the diagonal element in the stopping criteria provides a better in indication of the 
relative decrease in norm.
This could be implemented by changing line~\ref{alg2.line4} in Algorithm~\ref{alg:parlett} to
\[(c^p + r^p + |a_{ii}|^p) < 0.95 \times ( s + |a_{ii}|^p ).  \]
However, it is easier to absorb the addition of the diagonal element into the calculation of $c$ and $r$.
This is what we do in our proposed algorithm (Algorithm~\ref{alg:balance}). 
The only difference in Algorithm~\ref{alg:balance} and Algorithm~\ref{alg:parlett} is in line~\ref{alg3.line1}.
In Algorithm~\ref{alg:balance} we include the diagonal elements in the calculation of $c$ and $r$. 

\begin{algorithm}[htbp]
	\SetKwInOut{input}{Input}
	\SetKwInOut{output}{Output}
	\SetKwInOut{notes}{Notes}

	\Indm
	\input{ A matrix $A \in \R[n][n]$.  We will assume that $A$ is irrreducible.}
	\output{ A diagonal matrix $D$ and $A$ which is overwritten by $D^{-1} A D$. } 
	\notes{ $\beta$ is the radix base. 
		On output $A$ is nearly balanced in the $p$-norm. }
	\BlankLine
	\Indp
	$D \leftarrow I$ \\
	$converged \leftarrow 0$ \\
	\While{ $converged = 0$ } {
		$converged \leftarrow 1$ \\
		\For{ $i \leftarrow 1,\dots,n$ } {
			$\displaystyle c \leftarrow \| A(:,i) \|_p , \;\; 
			r \leftarrow \| A(i,:) \|_p $ \\ \label{alg3.line1}
			$s \leftarrow c^p + r^p , \;\;
			f \leftarrow 1 $ \\ 			
			\While{ $c < r/\beta$ } { \label{line2}
				$\displaystyle c \leftarrow c\beta, \;\;
				r \leftarrow r/\beta, \;\;
				f \leftarrow f \times \beta $ \\
			}
			\While{ $c \ge r\beta$ } {
				$\displaystyle c \leftarrow c/\beta, \;\;
				r \leftarrow r\beta, \;\;
				f \leftarrow f/\beta $ \\ \label{line3}
			} 
			\If{ $(c^p + r^p) < 0.95 \times s$ } { \label{line4}
				$ converged \leftarrow 0, \;\; 
				d_{ii} \leftarrow f \times d_{ii}$ \\
				$ A(:,i) \leftarrow f \times A(:,i), \;\;
				A(i,:) \leftarrow A(i,:)/f$ \\
			}
		}	
	}
	\caption{ Balancing (Proposed) } 
	\label{alg:balance}
\end{algorithm}
The 1-norm is desirable for efficiency and the typical 
choice of norm in the LAPACK library. 
However, we where unable to solve all of our test cases using the 1-norm.  
The 2-norm is able to balance the 
matrix while not destroying the backward error for all our test cases. 
Therefore, we propose using $p = 2$ in Algorithm~\ref{alg:balance}.

\section{Case Study}\label{sec:casestudy}
The balancing algorithm requires that the matrix be irreducible. The case study we present in this section shows the
potential danger if a matrix in nearly irreducible. 
Consider the matrix
\[
A=\left( \begin{array}{cccc} 1 & 1 & 0 & 0 \\ 0 & 2 & 1 & 0 \\ 0 & 0 & 3 & 1 \\ \epsilon & 0 & 0 & 4 \end{array} \right)
\]
where $0 \le \epsilon \ll 1$ \cite{LAPACK-forum-4270}.  As $\epsilon \rightarrow 0$ the eigenvalues tend towards the diagonal elements.  Although
these are not the exact eigenvalues (which are computed in Appendix~\ref{append:casestudy}), we will refer to the eigenvalues as the diagonal elements. 
The diagonal matrix 
\[
D = \alpha \left( \begin{array}{cccc} 
1 & 0 & 0 & 0 \\ 0 & \epsilon^{1/4} & 0 & 0 \\ 0 & 0 & \epsilon^{1/2} & 0 \\ 0 & 0 & 0 & \epsilon^{3/4}
 \end{array} \right),
\]
for any $\alpha > 0$, balances $A$ exactly for any $p$-norm.
The balanced matrix is
\[
\widetilde{A} = D^{-1}AD =\left( \begin{array}{cccc} 1 & \epsilon^{1/4} & 0 & 0 \\ 0 & 2 &  \epsilon^{1/4} & 0 \\ 
	0 & 0 & 3 &  \epsilon^{1/4} \\  \epsilon^{1/4} & 0 & 0 & 4 \end{array} \right).
\]

If $v$ is an eigenvector of $\widetilde{A}$, then $Dv$ is an eigenvector of $A$.   
The eigenvector associated with 4 is the most problematic.  The eigenspace of 4 for $A$ is 
$\textmd{span}( (1, 3, 6, 6)^T )$.  All of the components have
about the same magnitude and therefore contribute equally towards the backward error. 
The eigenspace of 4 for $\widetilde{A}$ approaches $\textmd{span}( (0, 0, 0, 1)^T )$ as $\epsilon \rightarrow 0$ and 
the components
begin to differ greatly in magnitude.  Error in the computation of the small elements are fairly irrelevant for 
backward error of $\widetilde{A}$.  But, when
converting back to the eigenvectors of $A$ the errors in the small components will be magnified and a poor 
backward error is expected. 
Computing the eigenvalue decomposition via Matlab's routine \texttt{eig} with $\epsilon = 10^{-32}$ the first component 
is zero. This results in
no accuracy in the eigenvector computation and backward error.   
Appendix~\ref{append:casestudy} provides detailed computations of the exact eigenvalues and eigenvectors of the 
balanced and unbalanced matrices.

This example can be generalized to more cases where one would expect the backward error to deteriorate.  
If a component of the eigenvector that is originally
of significance is made to be insignificant in the balanced matrix, then error in the computation is magnified 
when converting back to the original matrix.

Moreover, the balanced matrix fails to reduce the norm of the matrix by a significant amount since 
\[
	\frac{ \| \textmd{vec}( \widetilde{A} ) \|_1 } {  \| \textmd{vec}( A ) \|_1 } \rightarrow 10 / 13
\]
as $\epsilon \rightarrow 0$.  If we ignore the diagonal elements, then 
\[
	\frac{ \| \textmd{vec}( \widetilde{A}' ) \|_1 }{ \| \textmd{vec}( A' ) \|_1 } \rightarrow 0 ,
\] 
where $A'$ denotes the off diagonal entries of $A$.
Excluding the diagonal elements makes it appear as if balancing is reducing the norm of the matrix greatly.  But, when the 
diagonal elements are included the original matrix is already nearly balanced.

\section{Backward Error}\label{sec:backerror}
In this section $u$ is machine precision, $c(n)$ is a low order polynomial dependent only on the size of $A$, and $V_i$ is the
$i^{th}$ column of the matrix $V$.

Let $\widetilde{A} = D^{-1} A D$, where $D$ is a diagonal matrix and the elements of the diagonal of $D$ are restricted to being
exact powers of the radix base (typically 2).  Therefore, there is no error in computing $\widetilde{A}$.  
If a backward stable algorithm is used to compute the eigenvalue decomposition of $\widetilde{A}$, 
then $\widetilde{A}\widetilde{V} = \widetilde{V}\Lambda + E$, where $\| \widetilde{V}_i \|_2 = 1$ for all $i$
and $\| E \|_F \le c(n) u \| \widetilde{A} \|_F.$  The eigenvectors of $A$ are obtained by multiplying $\widetilde{V}$ by $D$, again there is 
no error in the computation $V = D\widetilde{V}$.  Scaling the columns of $V$ to have 2-norm of 1 gives
$A V_i = \lambda_{ii}V_i + D E_i / \| D \widetilde{V}_i \|_2$ and we have the following upper bound on the backward error:
\begin{align*}
	\| A V_i - \lambda_{ii}V_i \|_2 &= \frac{\| D E_i \|_2 }{\| D \widetilde{V}_i \|_2} \\
		&\le c(n) u \frac{\| D \|_2 }{\| D \widetilde{V_i} \|_2}\| \widetilde{A} \|_F  \\
		&\le c(n) u \frac{\| D \|_2 \| \widetilde{A} \|_F }{\sigma_{min}(D)} \\
		&= c(n) u \kappa(D) \| \widetilde{A} \|_F. 
\end{align*}
Therefore, 
\begin{equation}\label{eq:backerror}
	\frac{\| A V - V\Lambda \|_F}{\| A \|_F} \le c(n) u \frac{\kappa(D) \| \widetilde{A} \|_F}{\| A \|_F}. 
\end{equation}

If $\widetilde{V}_i$ is close to the singular vector associated with the minimum singular value then the bound 
will be tight.  This will be the case when
$\widetilde{A}$ is near a diagonal matrix.  On the other hand, if $\widetilde{A}$ is not near a diagonal matrix
then most likely $\| D \widetilde{V}_i \|_2 \approx \| D \|_2$ and the upper bound is no longer tight.  
For this case it would be 
better to expect the error to be close to $c(n) u \| \widetilde{A} \|_F$.

Ensuring that the quantity
\begin{equation}\label{bound}
u \frac{\kappa(D) \| \widetilde{A} \|_F }{ \| A \|_F}
\end{equation}
is small will guarantee that the backward error remains small.  The quantity is easy to compute
and therefore a reasonable stopping criteria for the balancing algorithm.  However, in many cases
the bound is too pessimistic and will tend to not allow balancing at all (see Section~\ref{sec:experiments}).

%
%
%
%
%
%

\section{Experiments}\label{sec:experiments}

We will consider four examples to justify balancing in the 2-norm and including the 
diagonal elements in the stopping criteria.  The first example is the case study 
discussed in Section~\ref{sec:casestudy} with $\epsilon = 10^{-32}$. 

The second example is a near triangular matrix (a generalization of the case study).  
Let $T$ be a triangular matrix with normally distributed nonzero entries.  Let $N$ be a dense random matrix with 
normally distributed entries.  Then $A = T + \epsilon N$, with $\epsilon = 10^{-30}$.  Balancing will tend towards
a diagonal matrix and the eigenvectors will contain entries with little significance.  Error in these components 
will be magnified when converting back to the eigenvectors of $A$.  

The third example is a Hessenberg matrix, that is, the Hessenberg form of a dense random matrix with normally 
distributed entries.  It is well known that balancing this matrix (using the LAPACK algorithm) 
will increase the eigenvalue condition number~\cite{watkins:2006:case}.

The finally example is a badly scaled matrix. Starting with a dense random matrix with normally 
distributed entries.  A diagonally similarity transformation with the logarithm of the 
diagonal entries evenly spaced between 0 and 10 gives us the initial matrix.  This is an example
where we want to balance the matrix.  This example is used to show that the new algorithm will
still balance matrices and improve accuracy of the eigenvalues.

To display our results we will display plot the quantities of interest versus iterations of the
balancing algorithm.  The black stars and squares on each line indicate the iteration in which the balancing 
algorithm converged.  Plotting this way allows use to see how the quantities vary during the 
balancing algorithm.  This is especially useful to determine if the error bound~\eqref{bound}
is a viable stopping criteria.

To check the accuracy of the eigenvalue we will examine the absolute condition number of an eigenvalue
which is given by the formula
\begin{equation}\label{condeig}
\kappa(\lambda,A) = \frac{ \| x \|_2 \| y \|_2 }{ | y^H x | },
\end{equation}
where $x$ and $y$ are the right and left eigenvectors associated with $\lambda$ \cite{kressner:2005:phd}.

\begin{figure}[!b]
	\centering
		\subfloat[Case study]{
			\resizebox{.4\textwidth}{!}{\includegraphics{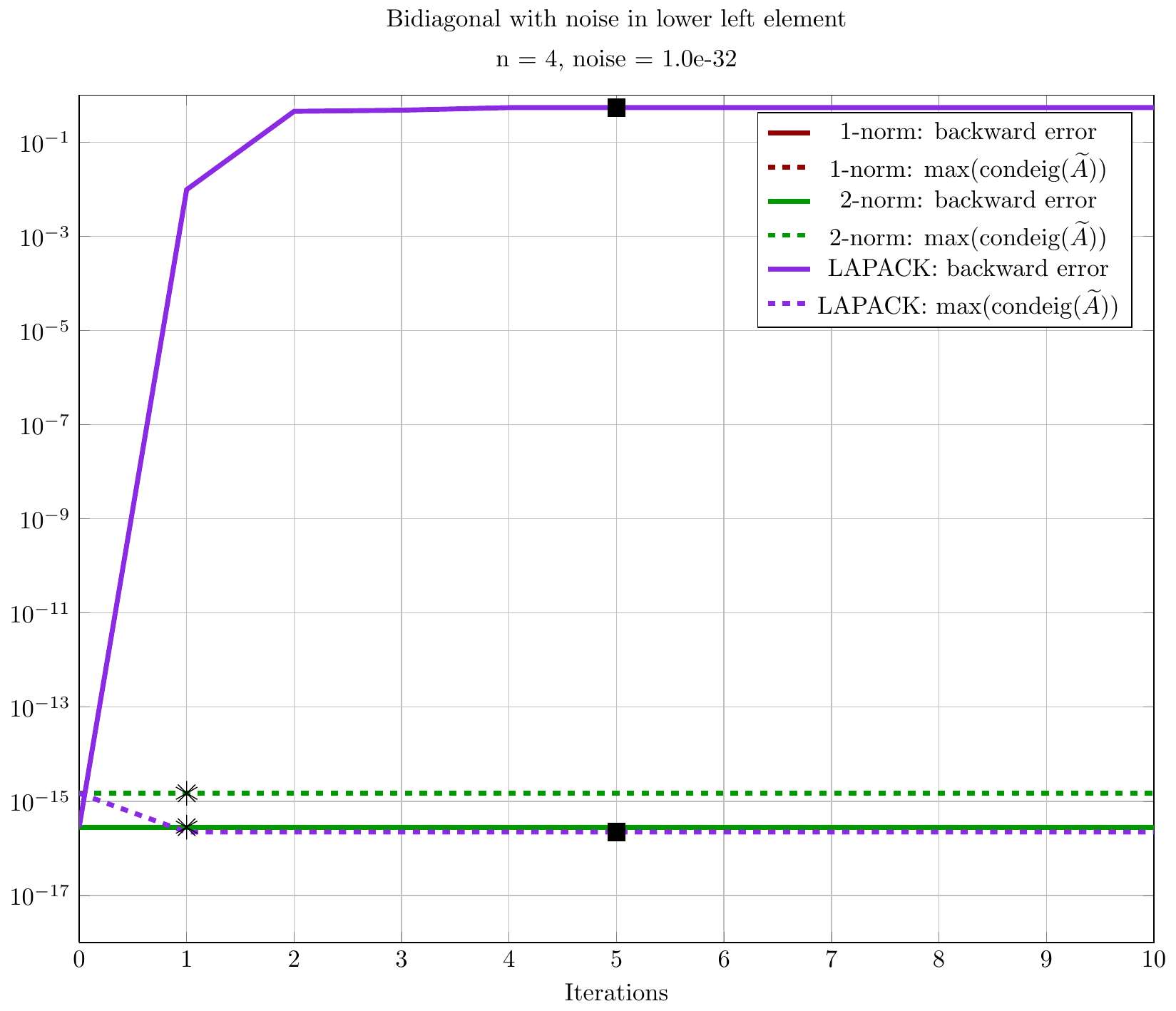}}
			\label{fig:lapackforum-1}
		}
		\hfil
		\subfloat[Near upper triangular]{
			\resizebox{.4\textwidth}{!}{\includegraphics{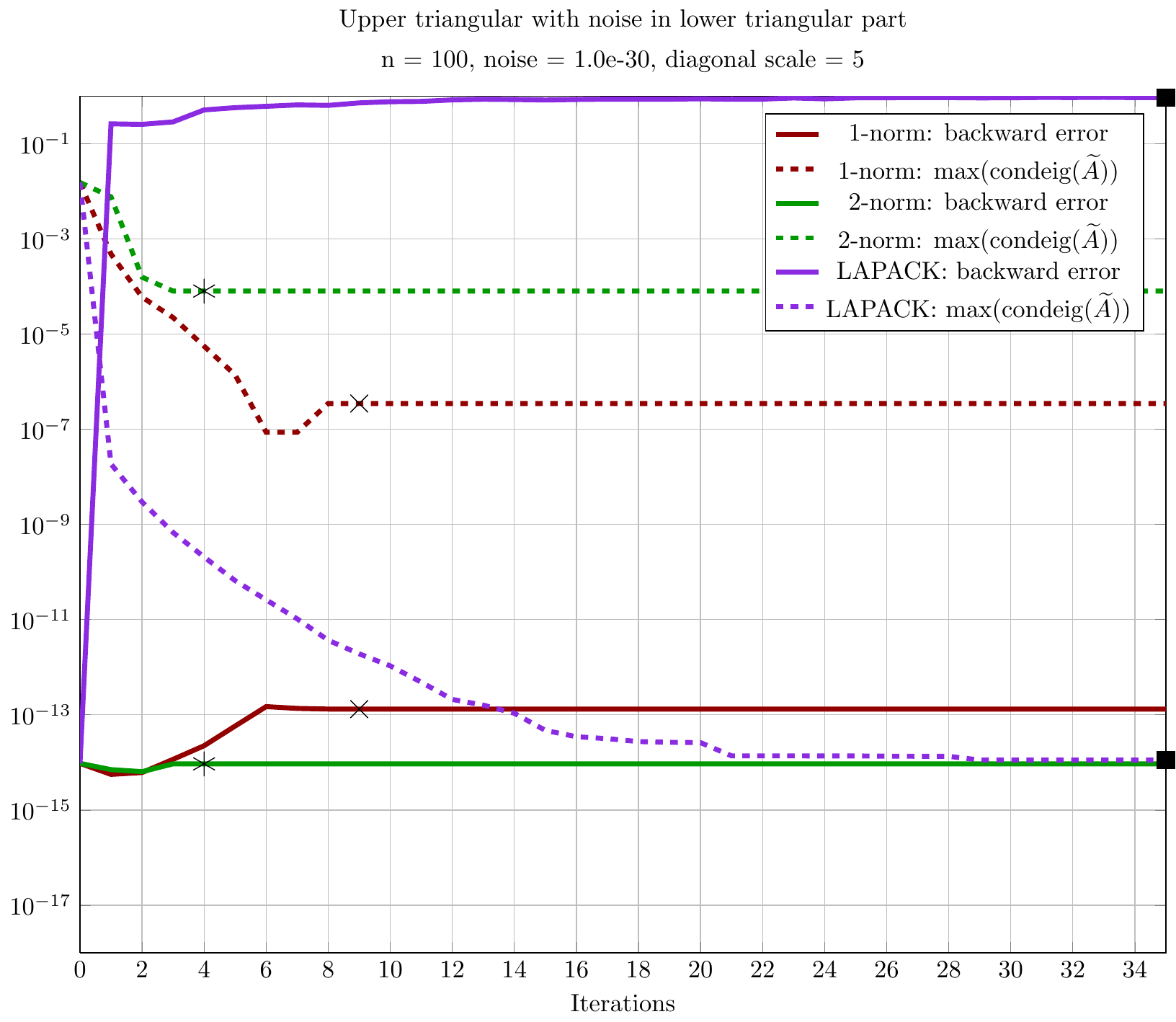}}
			\label{fig:uppertri-1}
		} \\
		\subfloat[Hessenberg form]{
			\resizebox{.4\textwidth}{!}{\includegraphics{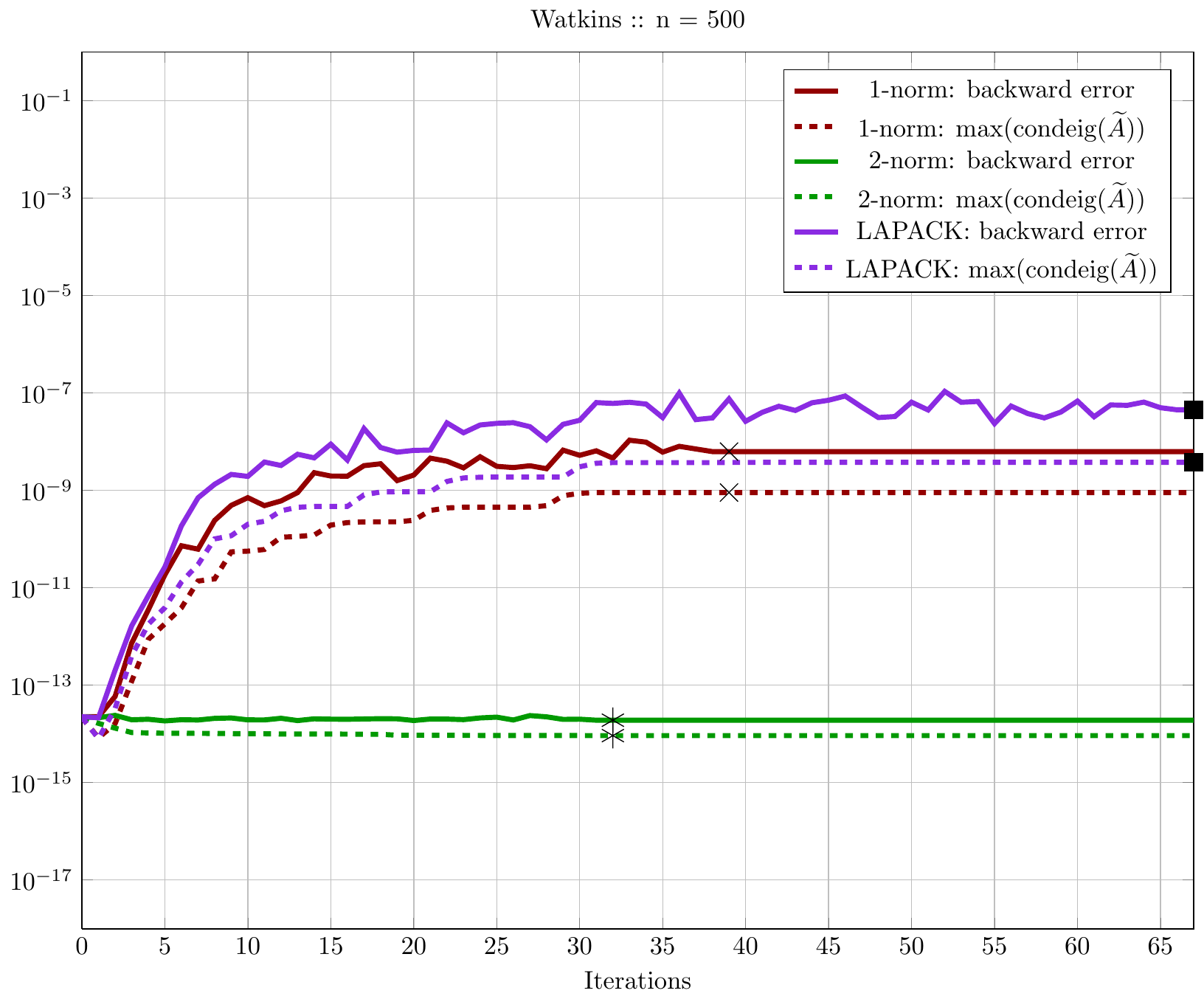}}
			\label{fig:watkins-1}
		}
		\hfil
		\subfloat[Badly scaled]{
			\resizebox{.4\textwidth}{!}{\includegraphics{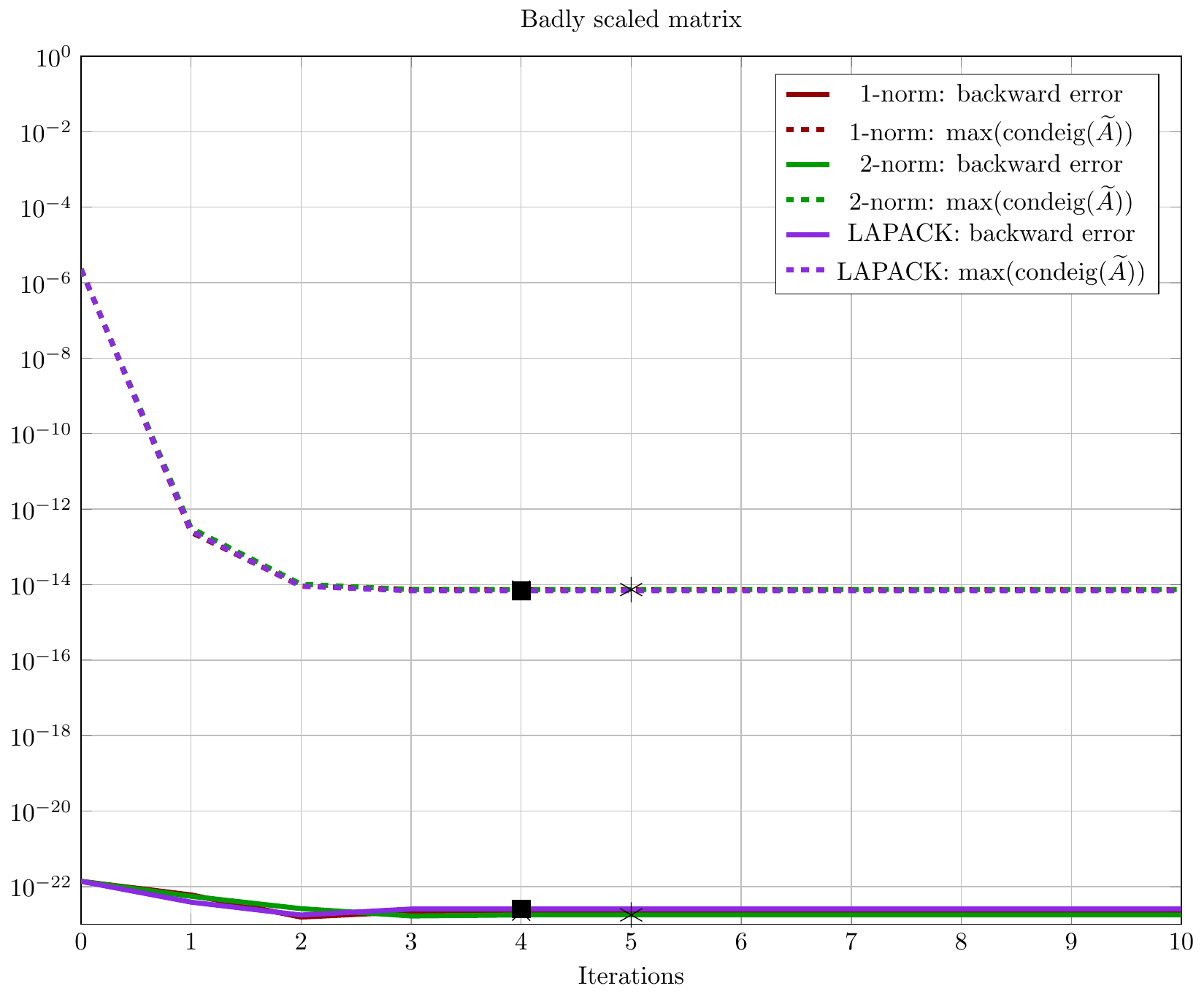}}
			\label{fig:badlyscaled-1}
		} 
		\caption{Backward error and absolute eigenvalue condition number versus iterations.}
		\label{fig:all-alg}
\end{figure}

Figure~\ref{fig:all-alg} plots the relative backward error, 
\begin{equation}\label{eq:relbackerror}
	\| A V - VD \|_2 / \| A \|_2,
\end{equation}
(solid lines) and machine precision times the
maximum absolute condition number of the eigenvalues, 
\[ u \cdot \max( \kappa(\lambda_{ii}, A)\; | \; 1 \le i \le n ), \] 
(dashed lines) versus iterations of the 
balancing algorithms. The dashed lines are an estimate of the forward error. 

In Figures~\ref{fig:lapackforum-1} and \ref{fig:uppertri-1} the LAPACK algorithm 
has no accuracy in the backward sense, while both the 1-norm and 2-norm algorithms maintain good
backward error.  
On the other hand the eigenvalue condition number is reduced much greater by the LAPACK algorithm in the
second example.  The LAPACK algorithm has a condition number 7 orders of magnitude less than the 
1-norm algorithm and 10 orders of magnitude less than the 2-norm algorithm.   
Clearly there is a trade-off to maintain desirable backward error.  If only the eigenvalues
are of interest, then the current balancing algorithm may still be desirable.  However, the special case of balancing
a random matrix already reduced to Hessenberg form would still be problematic and need to be avoided.  

In Figure~\ref{fig:watkins-1} both the backward error and the eigenvalue condition number worsen with the LAPACK and
1-norm algorithms.  The 2-norm algorithm does not allow for much scaling to occur and as a consequence
maintains a good backward error and small eigenvalue condition number.  This is the only
example for which the 1-norm is significantly worse than the 2-norm algorithm 
( 5 orders of magnitude difference).  This example is the only reason we choose to use the 2-norm over the 1-norm
algorithm.  

The final example (Figure~\ref{fig:badlyscaled-1}) shows all algorithms
have similar results:  maintain good backward error and reduce the eigenvalue condition number.
This indicates that the new algorithms would not prevent balancing when appropriate. 

\begin{figure}[!b]
	\centering
		\subfloat[Case study]{
			\resizebox{.4\textwidth}{!}{\includegraphics{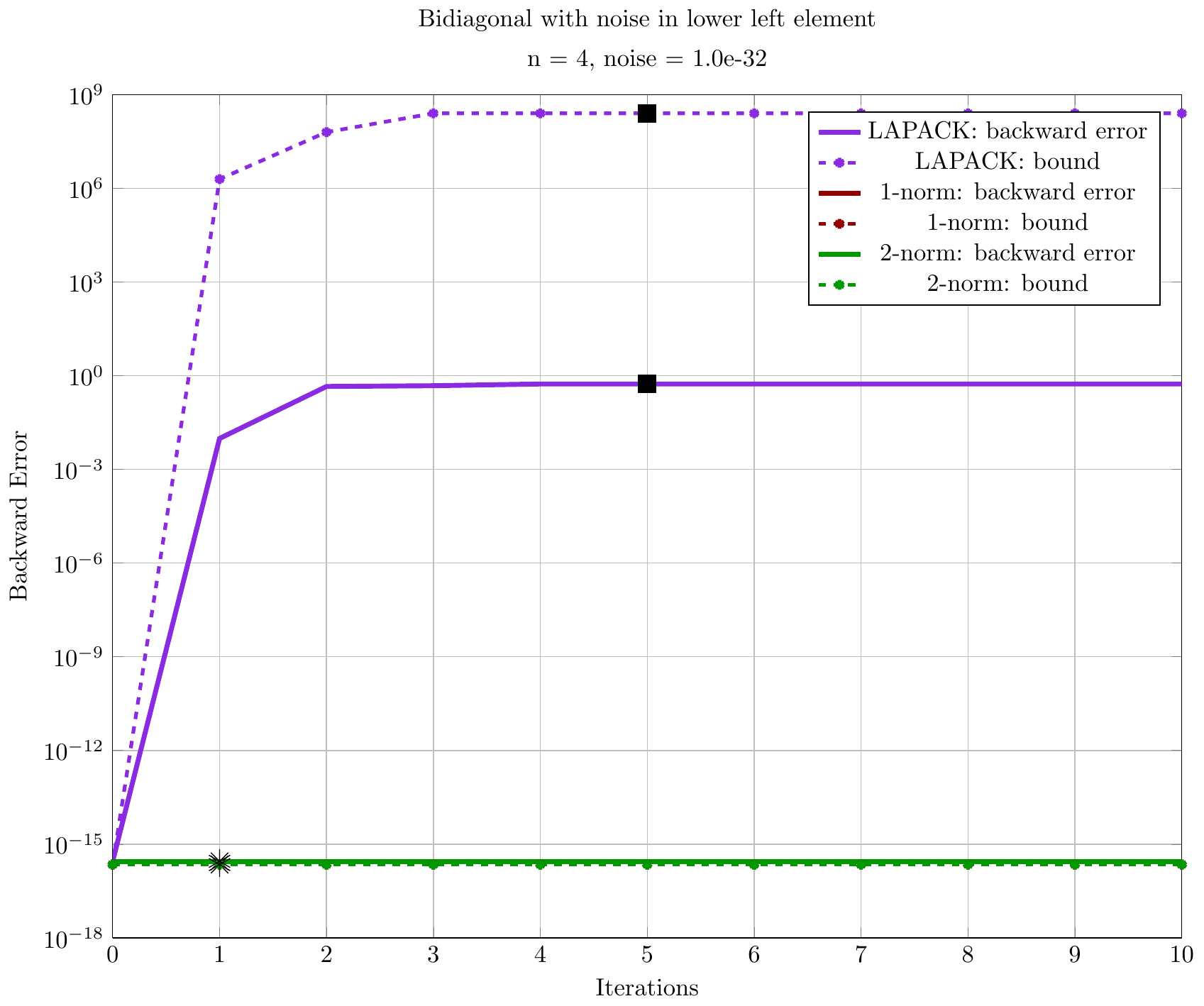}}
			\label{fig:lapackforum-2}
		}
		\hfil
		\subfloat[Near upper triangular]{
			\resizebox{.4\textwidth}{!}{\includegraphics{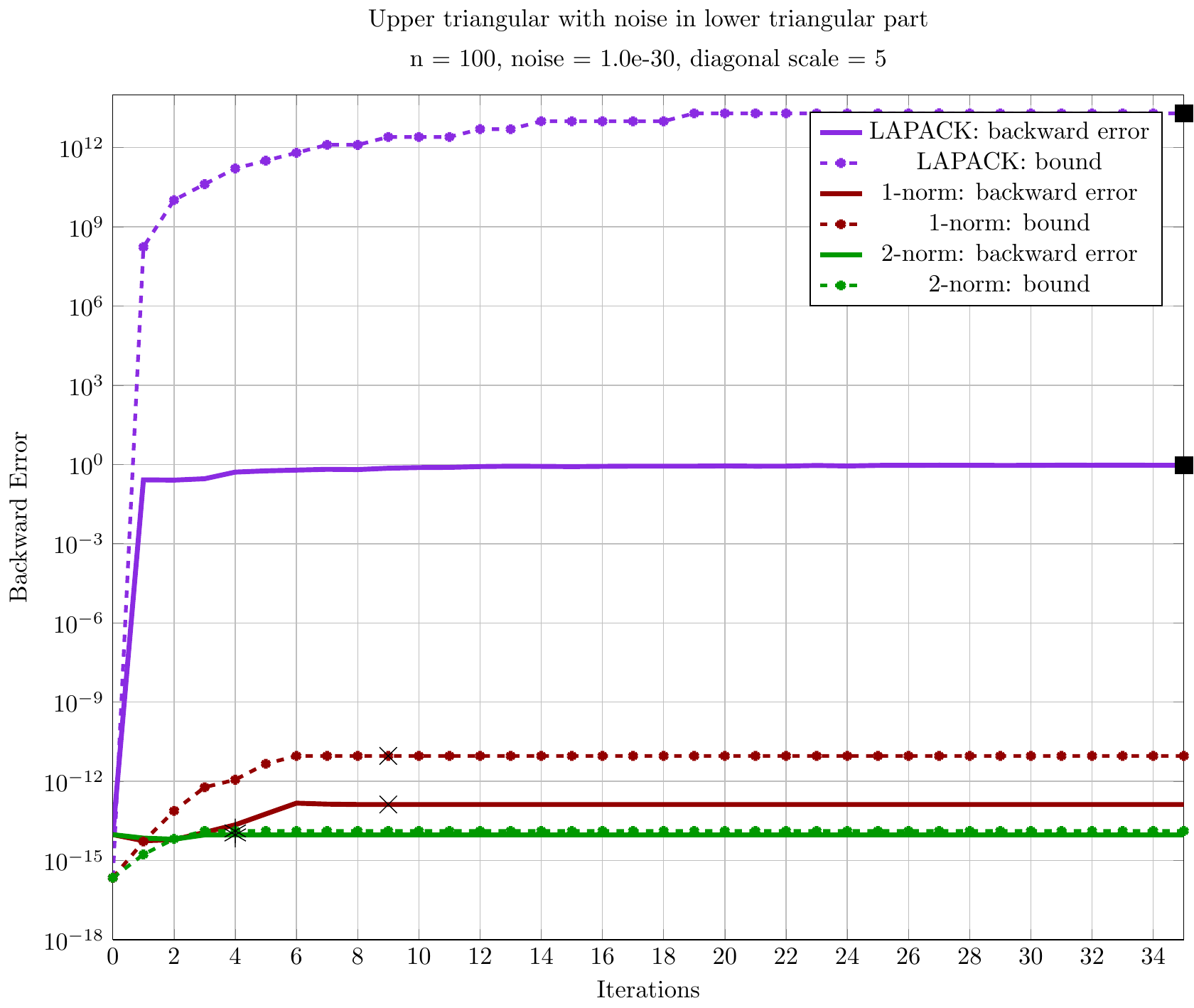}}
			\label{fig:uppertri-2}
		} \\
		\subfloat[Hessenberg form]{
			\resizebox{.4\textwidth}{!}{\includegraphics{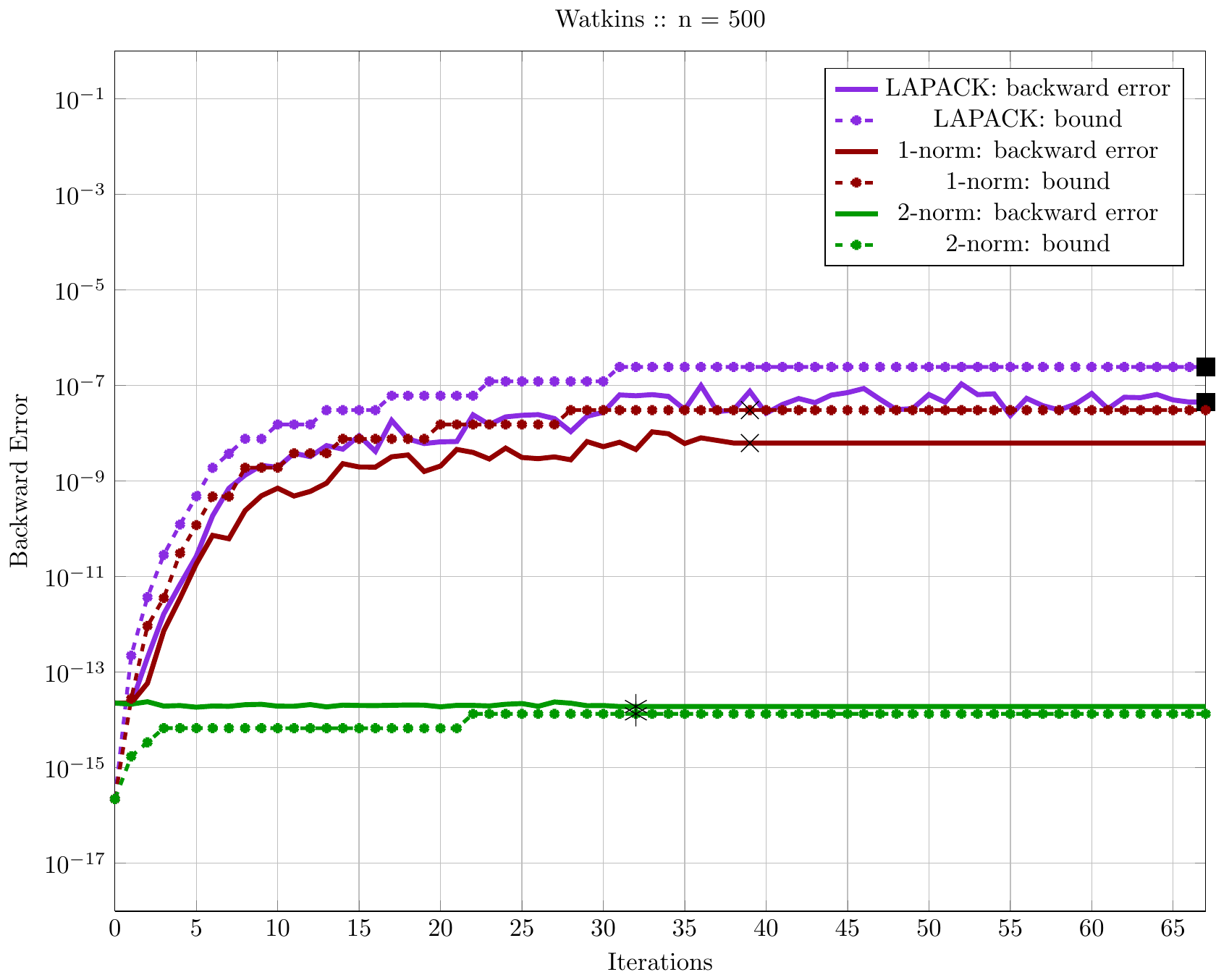}}
			\label{fig:watkins-2}
		}
		\hfil
		\subfloat[Badly scaled]{
			\resizebox{.4\textwidth}{!}{\includegraphics{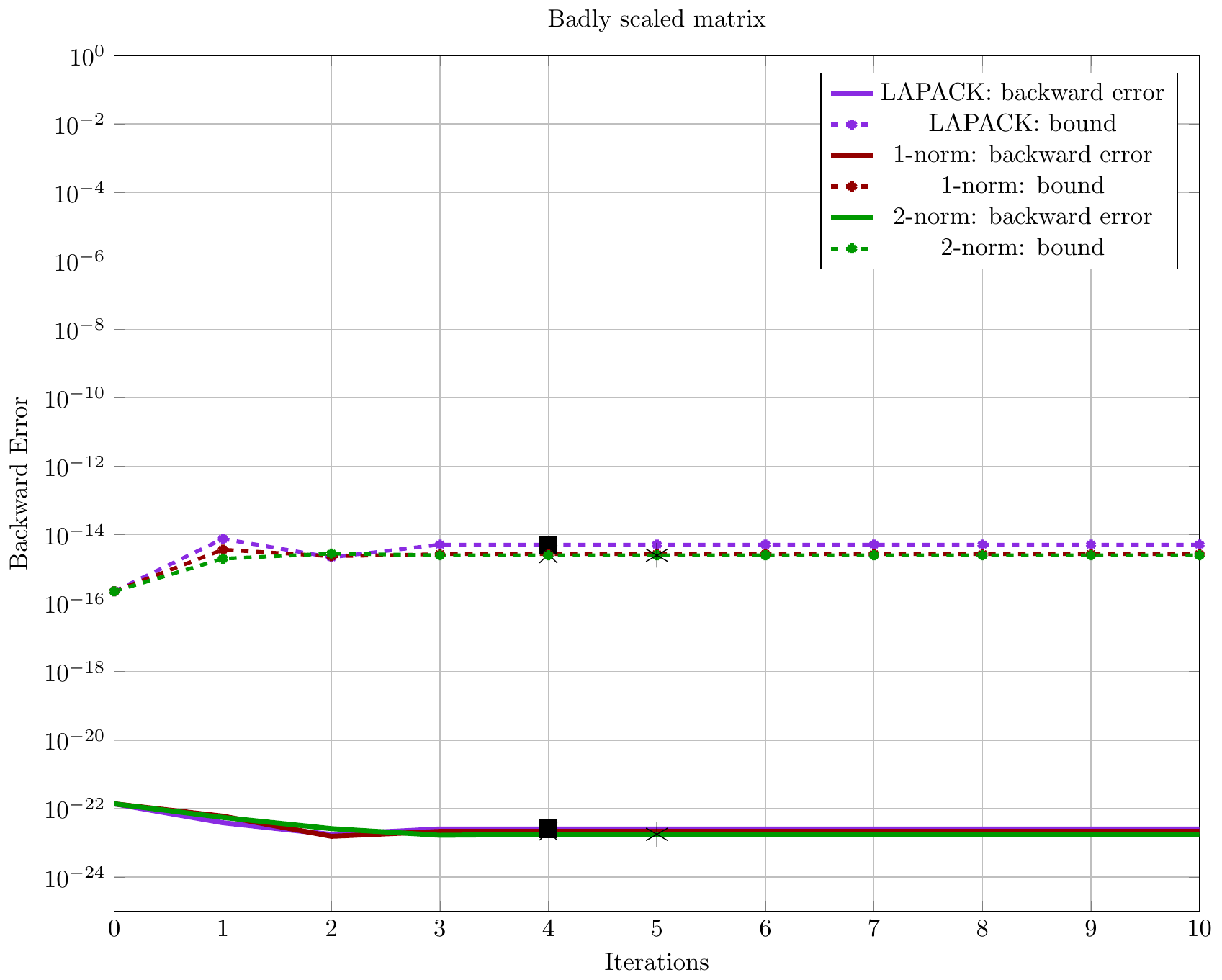}}
			\label{fig:badlyscaled-2}
		} 
		\caption{Backward error bounds. }
		\label{fig:bounds}
\end{figure}

Figure~\ref{fig:bounds} plots the relative backward error~\eqref{eq:relbackerror} (solid lines) 
and the backward error bound~\eqref{bound} (dashed lines).  
For the bound to be used as a stopping criteria, the bound should be
a fairly good estimate of the true error.  Otherwise, the bound would prevent balancing when balancing could have 
been used without hurting the backward error.  Figures~\ref{fig:lapackforum-2} and \ref{fig:uppertri-2} show the bound
would not be tight for the LAPACK algorithm and therefore not a good stopping criteria.
Figure~\ref{fig:badlyscaled-2} shows that the bound is not tight for all algorithms.  For this reason we do not 
consider using the bound as a stopping criteria.  

\begin{figure}[!t]
	\centering
		\subfloat[Case study]{
			\resizebox{.4\textwidth}{!}{\includegraphics{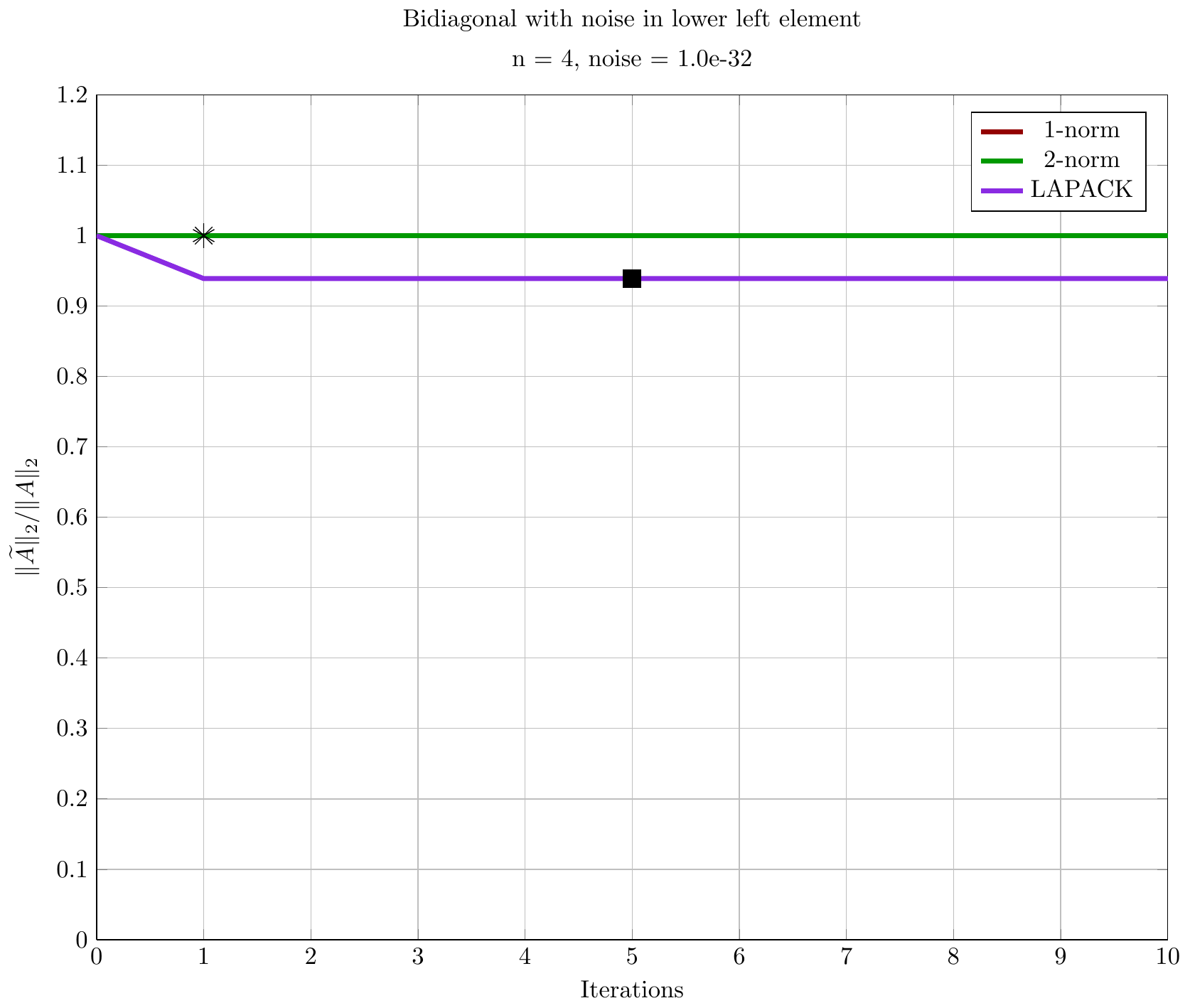}}
			\label{fig:lapackforum-3}
		}
		\hfil
		\subfloat[Near upper triangular]{
			\resizebox{.4\textwidth}{!}{\includegraphics{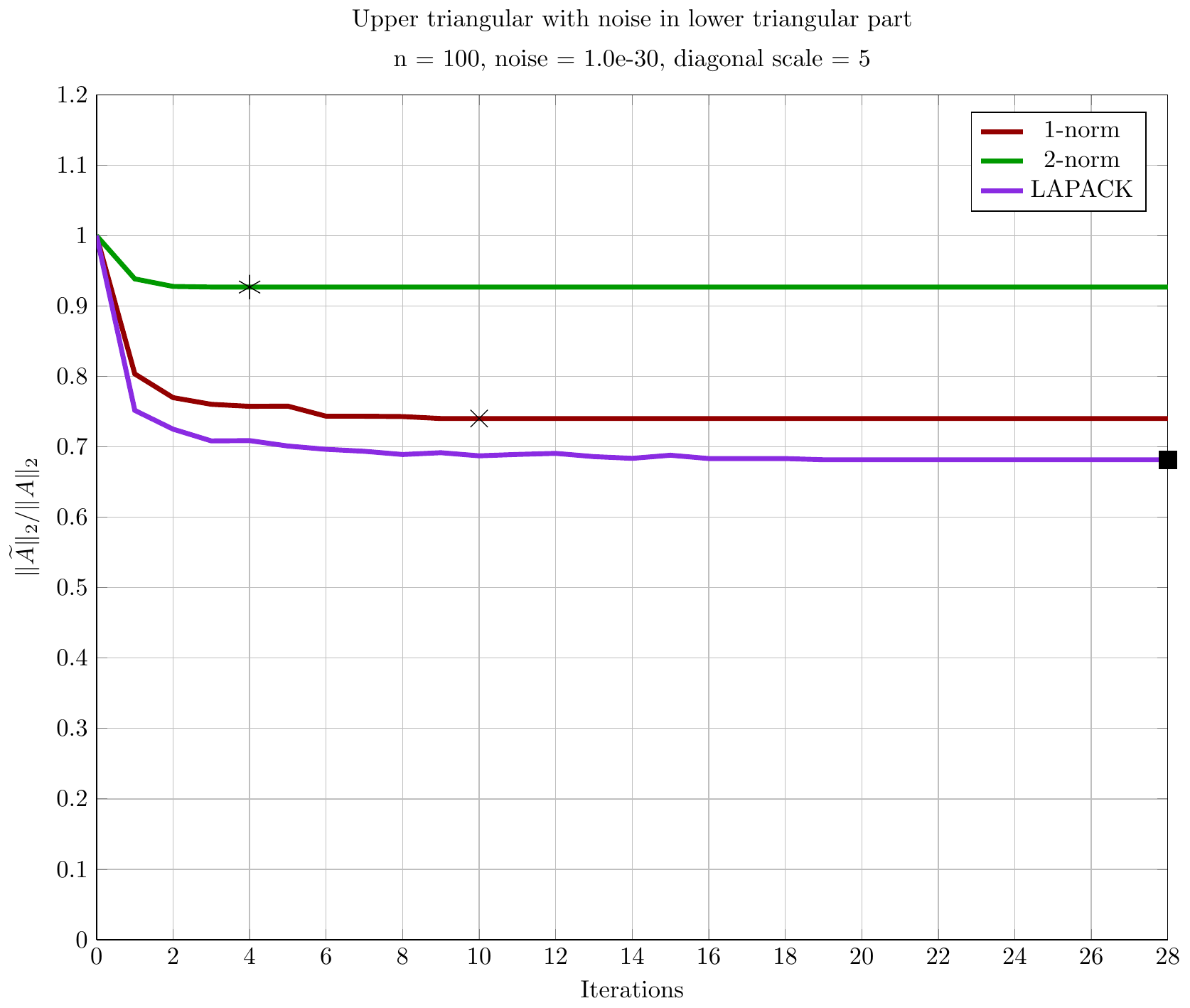}}
			\label{fig:uppertri-3}
		} \\
		\subfloat[Hessenberg form]{
			\resizebox{.4\textwidth}{!}{\includegraphics{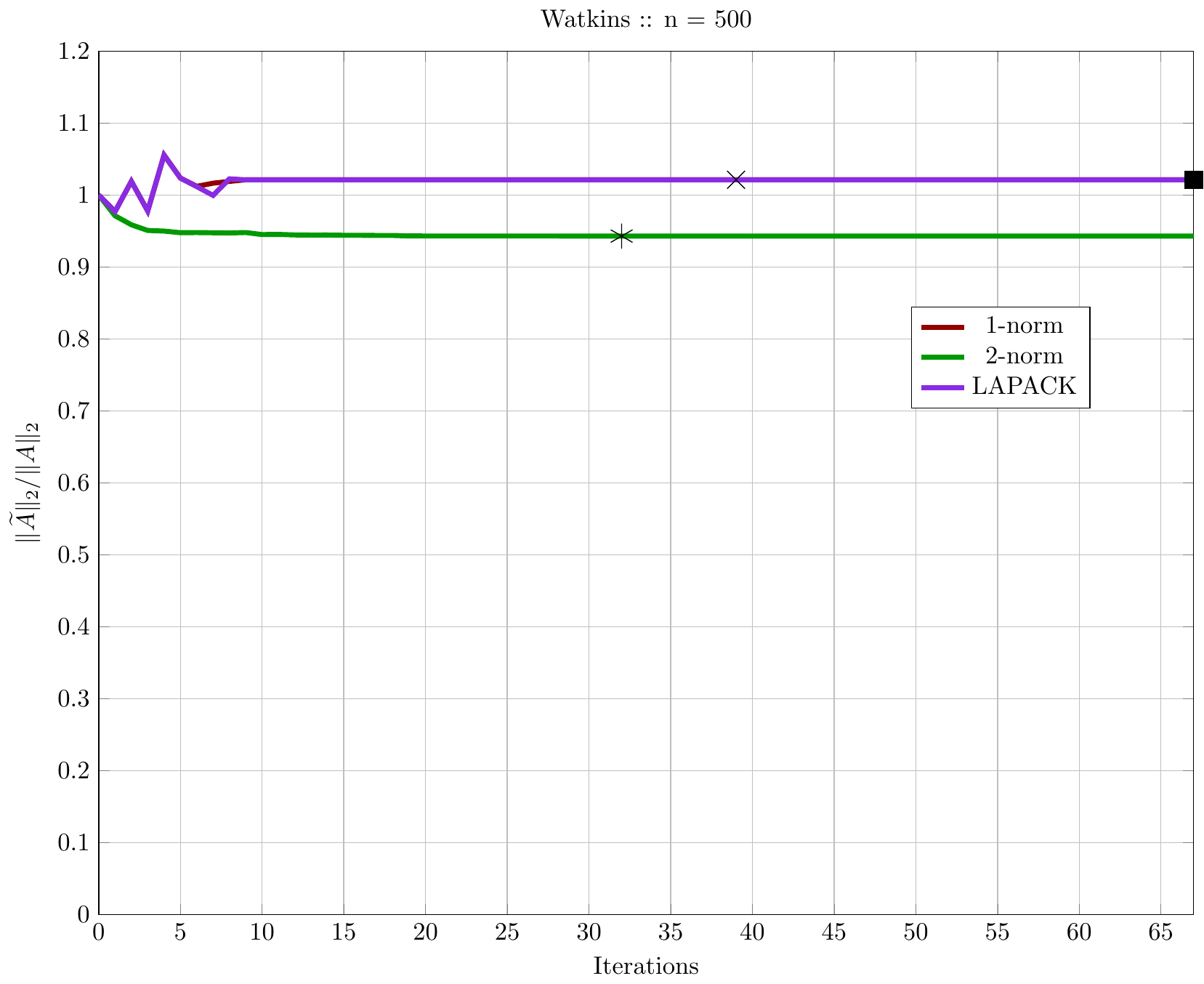}}
			\label{fig:watkins-3}
		}
		\hfil
		\subfloat[Badly scaled (log scale)]{
			\resizebox{.4\textwidth}{!}{\includegraphics{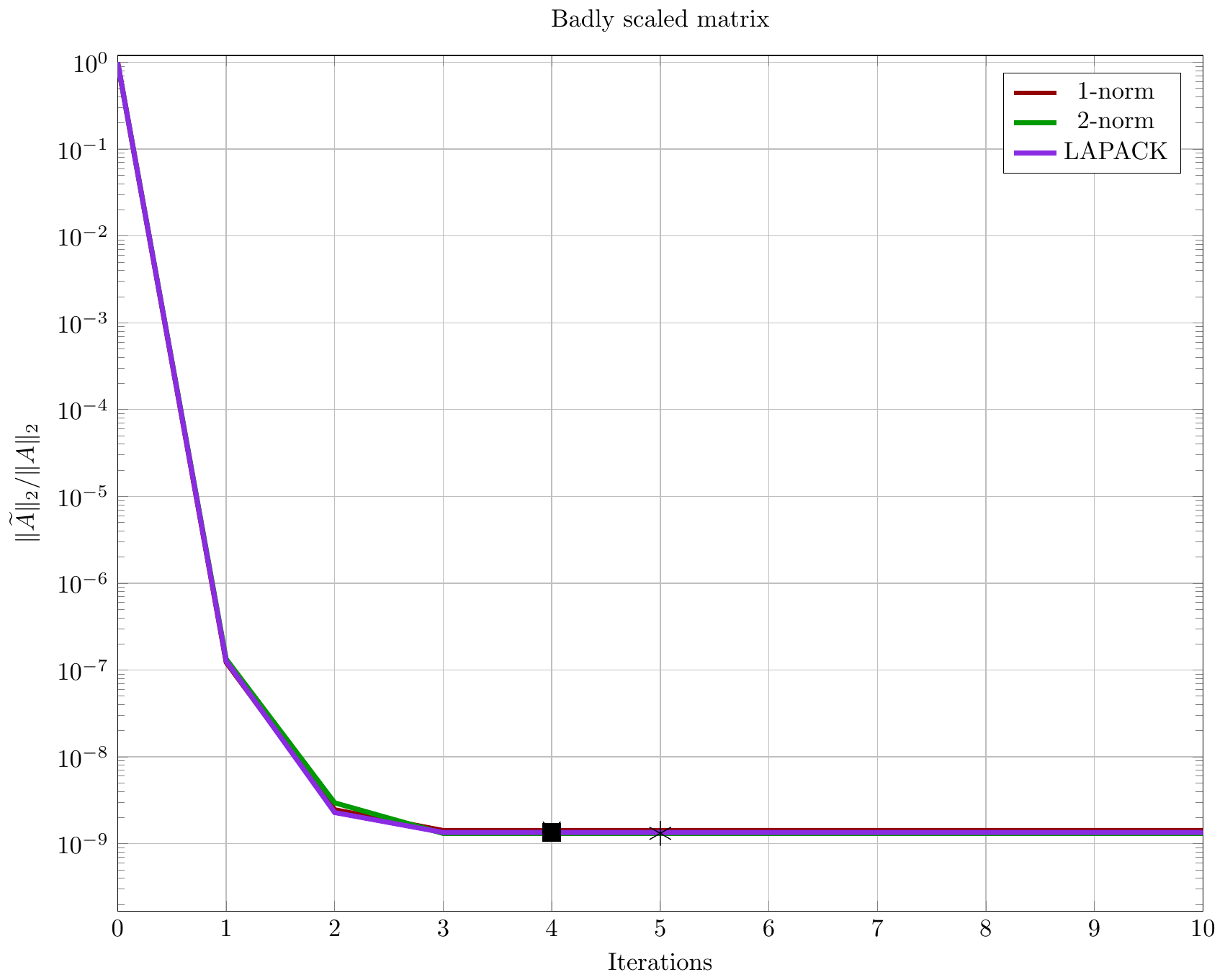}}
			\label{fig:badlyscaled-3}
		} 
		\caption{Ratio of balanced matrix norm and original matrix norm.}
		\label{fig:reducenorm}
\end{figure}

Figure~\ref{fig:reducenorm} shows the ratio of the norm of the balanced matrix and the 
norm of the original matrix versus iterations of the balancing algorithm.  
In Figures~\ref{fig:lapackforum-3},~\ref{fig:uppertri-3},~\ref{fig:watkins-3} the norm is not decreased by much in any of the
algorithms.  The maximum decrease is only 70\% for the near triangular matrix (Figure~\ref{fig:uppertri-3}).
In Figure~\ref{fig:badlyscaled-3} all algorithms successfully reduce the norm of the matrix by 9 
orders of magnitude.
In all three examples where the LAPACK algorithm deteriorates the backward error, the norm
is not successfully reduced by a significant quantity.  But, in our example where balancing is most 
beneficial the norm is reduced significantly.  

\section{Matrix Exponential}

Balancing could also be used for the computation of the matrix exponential, since
$e^{XAX^{-1}} = Xe^{A}X^{-1}$.  If the exponential of the balanced matrix can be
computed more accurately, then balancing may be justified.  There are numerous ways to
compute the matrix exponential.  The methods are discussed in the paper by  Moler and Van Loan,
``Nineteen Dubious Ways to Compute the Exponential of a Matrix''~\cite{MolerVanLoan:1978} and later
updated in 2003~\cite{MolerVanLoan:2003}.  The scaling and squaring method is consider one of the 
best options.  Currently, MATLAB's function \texttt{expm} implements this method 
using a variant by Higham~\cite{Higham:JMA:2005}.  The scaling and squaring 
method relies on the Pad{\'e} approximation of the matrix exponential, 
which is more accurate when $\| A \|$ is small. Hence, balancing
would be beneficial. 
 
However,  the matrix exponential possess the special property: $e^A = ( e^{A/s} )^s$. 
Scaling the matrix by $s$ successfully reduces the norm of $A$ to an acceptable level for an
accurate computation of $e^{A/s}$ by the Pad{\'e} approximate.
When $s$ is a power of 2 then repeat squaring is used to obtain the final result.   
Since the scaling of $A$ successfully achieves the same goal as balancing, balancing is not required
as a preprocessing step.  

\section{Conclusion}

We proposed a simple change to the LAPACK balancing algorithm \texttt{GEBAL}.  Including the diagonal elements
in the balancing algorithm provides a better measure of significant decrease in the matrix norm than the current 
version.  By including the diagonal elements and switching to the 2-norm we are able to avoid deteriorating the 
backward error for our test cases. 

The proposed change is to protect the backward error of the computation.  In some case the new algorithm 
will prevent from decreasing the maximum eigenvalue condition number as much as the current algorithm.
If only the eigenvalues are of interest, then it may be desirable to uses the current version.  
However, for the case where $A$ is dense and poorly scaled,
the new algorithm will still balance the matrix and improve the eigenvalue condition number.
If accurate eigenvectors are desired, then one should consider not balancing the matrix.

\appendix
\section{Detailed Case Study}\label{append:casestudy}
Consider the matrix
\[
A=\left( \begin{array}{cccc} 1 & 1 & 0 & 0 \\ 0 & 2 & 1 & 0 \\ 0 & 0 & 3 & 1 \\ \epsilon & 0 & 0 & 4 \end{array} \right)
\]
where $0 \le \epsilon \ll 1$.
The eigenvalues of $A$ are
\[
\Lambda=\left(
\begin{array}{cccc}
\frac{1}{2}\left(5-\sqrt{5+4\rope}\right)&0&0&0\\
0&\frac{1}{2}\left(5-\sqrt{5-4\rope}\right) &0&0\\ 
0&0& \frac{1}{2}\left(5+\sqrt{5-4\rope}\right)&0\\
0&0&0& \frac{1}{2}\left(5+\sqrt{5+4\rope}\right)\\
\end{array}
\right)\]
with corresponding eigenvectors
\[ 
V = \left(
\begin{array}{cccc}
\frac{2}{1+\rope}                         & \frac{2}{\left(3-\sqrt{5-4\rope}\right)}  & \frac{2}{\left(3+\sqrt{5-4\rope}\right)}  & \frac{2}{\left(3+\sqrt{5+4\rope}\right)} \\
\frac{3-\sqrt{5+4\rope}}{1+\rope}         & 1                                         & 1                                         & 1 \\
\frac{4+2\rope-2\sqrt{5+4\rope}}{1+\rope} & \frac{1}{2}\left(1-\sqrt{5-4\rope}\right) & \frac{1}{2}\left(1+\sqrt{5-4\rope}\right) & \frac{1}{2}\left(1+\sqrt{5+4\rope}\right) \\
3-\sqrt{5+4\rope}                         & 1-\rope                                   & 1-\rope                                   & 1+\rope
\end{array}
\right).
\]
Note that as $\epsilon \rightarrow 0$, we have 
\[
\Lambda \rightarrow \left(
\begin{array}{cccc}
 1&0&0&0 \\ 
 0&2&0&0 \\
 0&0&3&0 \\
 0&0&0&4 \\
\end{array}
\right)
\]
and (with some rescaling)
\[
V \rightarrow \left(
\begin{array}{cccc}
1 & 1 & 1 & 1 \\
0 & 1 & 1 & 3 \\
0 & 0 & 2 & 6 \\
0 & 0 & 0 & 6
\end{array}
\right).
\]

To balance $A$, we must find an invertible diagonal matrix 
\[
D = \left( \begin{array}{cccc} 
\alpha & 0 & 0 & 0 \\ 0 & \beta & 0 & 0 \\ 0 & 0 & \gamma & 0 \\ 0 & 0 & 0 & \delta
 \end{array} \right)
\]
such that $\widetilde{A}=D^{-1}A D$ has equal row and column norms for each diagonal element, where

\[D^{-1}A D = 
\left(
\begin{array}{cccc}
 1 & \beta/\alpha & 0 & 0 \\
 0 & 2 & \gamma/\beta & 0 \\
 0 & 0 & 3 & \delta/\gamma \\
  \epsilon \alpha/\delta & 0 & 0 & 4 \\
\end{array}
\right).\]
If we require that the diagonal elements of $D$ are all positive, the balancing condition can be written as 
\begin{equation} \label{bal}
\frac{\delta}{\gamma}=\frac{\gamma}{\beta}=\frac{\beta}{\alpha} = \epsilon \frac{\alpha}{\delta}.
\end{equation}
The balancing condition (\ref{bal}) can be met by parameterizing $D$ in terms of $\alpha$ where
\begin{eqnarray*}
\beta &=& \epsilon^{1/4} \alpha \\
\gamma &=& \epsilon^{1/2} \alpha \\
\delta &=& \epsilon^{3/4} \alpha
\end{eqnarray*}
and then the balancing matrix $D$ is
\[
D = \alpha \left( \begin{array}{cccc} 
1 & 0 & 0 & 0 \\ 0 & \epsilon^{1/4} & 0 & 0 \\ 0 & 0 & \epsilon^{1/2} & 0 \\ 0 & 0 & 0 & \epsilon^{3/4}
 \end{array} \right).
\]
The balanced matrix $\widetilde{A}$ then
\[
\widetilde{A} = D^{-1}AD =\left( \begin{array}{cccc} 1 & \epsilon^{1/4} & 0 & 0 \\ 0 & 2 &  \epsilon^{1/4} & 0 \\ 0 & 0 & 3 &  \epsilon^{1/4} \\  \epsilon^{1/4} & 0 & 0 & 4 \end{array} \right).
\]

Computing the eigenvectors of $\widetilde{A}$ directly, we have
\[
\widetilde{V}= \left( \begin{array}{cccc}
\frac{-3-\sqrt{5+4\rope}}{2}                                 & -\frac{\epsilon^{1/4}\left({3+\sqrt{5-4 \rope}}\right)}{2} & \frac{\epsilon^{1/2}\left( -3+\sqrt{5-4 \rope}\right)}{2} & \frac{2 \epsilon^{3/4}}{\left(1+\rope\right) \left(3+\sqrt{5+4\rope}\right)}\\
\frac{\epsilon^{3/4}}{1+\rope}                               & -1-\rope                                                   & -\epsilon^{1/4}\rope                                      & \frac{\epsilon^{1/2}}{1+\rope}\\
-\frac{\epsilon^{1/2}}{2} \frac{-1+\sqrt{5+4\rope}}{1+\rope} & -\frac{2\epsilon^{3/4}}{1+\sqrt{5-4\rope}}                 & -\frac{\rope \left( 1+ \sqrt{5-4\rope} \right)}{2}        & \frac{2 \epsilon^{1/4} \left(2+\rope+\sqrt{5+4 \rope}\right)}{\left(1+\rope\right) \left(3+\sqrt{5+4 \rope}\right)}\\
\epsilon^{1/4}                                               & \epsilon^{1/2}                                             & \epsilon^{3/4}                                            & 1\\
\end{array} \right).
\]
Note that as $\epsilon \rightarrow 0$, the balanced matrix $\widetilde{A}$ tends towards a diagonal matrix:
\[
\widetilde{A} \rightarrow \left( \begin{array}{cccc} 1 & 0 & 0 & 0 \\ 0 & 2 &  0 & 0 \\ 0 & 0 & 3 & 0 \\  0 & 0 & 0 & 4 \end{array} \right)
\]
with the eigenvectors of the balanced matrix approaching the identity matrix.

\nocite{lemeire:1976:BIT}
\nocite{watkins:2007:MEP}
\nocite{handbookautocomp:1971:springer}

\bibliographystyle{plain}
\bibliography{biblio_balancing}

\begin{thebibliography}{10}

\bibitem{templatesSolutions:2000:siam}
Zhaojun Bai, James Demmel, Jack Dongarra, Axel Ruhe, and Henk van~der Vorst,
  editors.
\newblock {\em Templates for the Solution of Algebraic Eigenvalue Problems}.
\newblock SIAM, Philadelphia, 2000.

\bibitem{chen:1998:masters}
Tzu-Yi Chen.
\newblock Balancing sparse matrices for computing eigenvalues.
\newblock Master's thesis, University of California at Berkeley, 1998.

\bibitem{chen:2001:phd}
Tzu-Yi Chen.
\newblock {\em Preconditioning Sparse Matrices for Computing Eigenvalues and
  Solving Linear Systems of Equations}.
\newblock PhD thesis, University of California at Berkeley, 2001.

\bibitem{chen:2000:LAA}
Tzu-Yi Chen and James~W Demmel.
\newblock Balancing sparse matrices for computing eigenvalues.
\newblock {\em Linear Algebra and Its Applications}, 309(1):261--287, 2000.

\bibitem{Higham:JMA:2005}
Nicholas~J Higham.
\newblock The scaling and squaring method for the matrix exponential revisited.
\newblock {\em SIAM Journal on Matrix Analysis and Applications},
  26(4):1179--1193, 2005.

\bibitem{kressner:2005:phd}
D.~Kressner.
\newblock {\em Numerical methods for general and structured eigenvalue
  problems}, volume~46 of {\em Lecture Notes in Computational Science and
  Engineering Series}.
\newblock Springer-Verlag Berlin and Heidelberg GmbH \& Company KG, 2005.

\bibitem{lemeire:1976:BIT}
Frans Lemeire.
\newblock Equilibration of matrices to optimize backward numerical stability.
\newblock {\em BIT Numerical Mathematics}, 16(2):143--145, 1976.

\bibitem{MolerVanLoan:1978}
Cleve Moler and Charles~Van Loan.
\newblock Nineteen dubious ways to compute the exponential of a matrix.
\newblock {\em SIAM Review}, 20(4):pp. 801--836, 1978.

\bibitem{MolerVanLoan:2003}
Cleve Moler and Charles~Van Loan.
\newblock Nineteen dubious ways to compute the exponential of a matrix,
  twenty-five years later.
\newblock {\em SIAM Review}, 45(1):pp. 3--49, 2003.

\bibitem{LAPACK-forum-4270}
Ron Morgan.
\newblock {http://icl.cs.utk.edu/lapack-forum/viewtopic.php?f=13\&t=4270},
  2013.

\bibitem{Osborne:1960:PM:321043.321048}
E.~E. Osborne.
\newblock On pre-conditioning of matrices.
\newblock {\em J. ACM}, 7(4):338--345, October 1960.

\bibitem{parlett:1969:NM}
Beresford~N Parlett and Christian Reinsch.
\newblock Balancing a matrix for calculation of eigenvalues and eigenvectors.
\newblock {\em Numerische Mathematik}, 13(4):293--304, 1969.

\bibitem{watkins:2006:case}
David~S. Watkins.
\newblock A case where balancing is harmful.
\newblock {\em Electron. Trans. Numer. Anal}, 23:1--4, 2006.

\bibitem{watkins:2007:MEP}
David~S. Watkins.
\newblock {\em The Matrix Eigenvalue Problem: GR and Krylov Subspace Methods}.
\newblock Society for Industrial and Applied Mathematics, Philadelphia, PA,
  USA, 1 edition, 2007.

\bibitem{handbookautocomp:1971:springer}
J.~H. Wilkinson and C.~Reinsch, editors.
\newblock {\em Handbook for Automatic Computation, Volume II, Linear Algebra}.
\newblock Springer-Verlag, New York, 1971.

\end{thebibliography}

\end{document}